\documentclass[11pt]{amsart}
\usepackage[section]{placeins}   

\usepackage{amsmath}
\usepackage{amsthm}
\usepackage{amssymb}
\usepackage[T1]{fontenc}
\usepackage{ae,aecompl}
\usepackage[all, knot, poly]{xy}
\input xy
\xyoption{all}

\newtheorem{theorem}{Theorem}[section] 
\newtheorem*{theorem*}{Theorem} 
 
\newtheorem{lemma}[theorem]{Lemma} 
\newtheorem{conjecture}[theorem]{Conjecture} 
\newtheorem*{conj*}{Conjecture} 
\newtheorem{prop}[theorem]{Proposition} 
\newtheorem{definition}[theorem]{Definition} 
 
\newtheorem{exa}[theorem]{Example}

\newtheorem*{ques*}{Question}

\def\inv{^{-1}}

\def\R{\mathbb{R}}

\def\F{\mathbb{F}}

\def\M{\mathcal{M}}

\def\comp{\circ}
\def\set#1{\left\{ #1\right\}}
\def\ie{\textit{i.e.\ }}

\def\h{\frac{1}{2}}
\def\inclusion#1#2{\xymatrix{ #1\ar@{^{(}->}[r] & #2}}
\def\incl{\ar@{^{(}->}}
\def\maps{\ar@{|->}}
\def\monic{\ar@{>->}}

\def\bd{\begin{defn}}
\def\ed{\end{defn}}
\def\ec{\end{cor}}
\def\bc{\begin{cor}}
\def\bl{\begin{lem}}
\def\el{\end{lem}}
\def\bt{\begin{thm}}
\def\et{\end{thm}}
\def\bex{\begin{exa}}
\def\eex{\end{exa}}
\def\bp{\begin{prop}}
\def\ep{\end{prop}}
\def\ben{\begin{enumerate}}
\def\een{\end{enumerate}}
\def\bi{\begin{itemize}}
\def\ei{\end{itemize}}
\def\be{\begin{equation}}
\def\bpm{\begin{pmatrix}}
\def\epm{\end{pmatrix}}
\def\ee{\end{equation}}

\def\cofun#1#2{\xymatrix{ #1 \ar@{ {~}{>}}[r] &  #2}  }
\def\confun#1#2{\xymatrix{ #1&  #2 \ar@{ {~}{>}}[l] }  }

\def\nullset{\emptyset}

\newcommand{\catname}[1]{\mbox{\sffamily\upshape {#1}}}
\newcommand{\Link}{\catname{Link}}
\newcommand{\Vect}{\catname{Vect}}
\newcommand{\Mod}{\catname{Mod}}

\usepackage{url, enumitem, amscd}
\usepackage{geometry}
\usepackage{graphicx}
\graphicspath{{Pictures/}}
\usepackage[usenames,dvipsnames]{color}
\usepackage[T1]{fontenc}
\usepackage{lmodern}

\def\K{{\mathcal K}}
\def\M{{\mathcal M}}

\def\D{{\mathcal D}}

\def\MM{{\mathfrak M}}

\def\F{\mathbb{F}}

\def\og{\overline{g}}

\DeclareMathOperator*{\HFL}{HFL}

\DeclareMathOperator*{\Kh}{Kh}

\title[Movie Moves for Knotted Surfaces with Markings]{Movie Moves for Knotted Surfaces with Markings}
\author[Matthew  Graham]{Matthew  Graham}
\email{mdgraham@neiu.edu}
\address{Department of Mathematics, 
Northeastern Illinois  University
5500 St Louis Ave, Chicago, IL 60625
}

\begin{document}

\maketitle


\begin{abstract}
 We present a marked analogue of Carter and Saito's  movie theorem.  Our  definition of marking was chosen to coincide with the markings that arise in link Floer homology.  In order to deal with complications arising from certain isotopies, we define three equivalences for marked surfaces and work over an equivalence class of marked surfaces when proving our generalization of Carter and Saito's movie theorem.
\end{abstract}

\section{Introduction}

Our initial motivation to study marked surfaces came from link Floer homology ($\HFL$).  
Link Floer homology, is a powerful link invariant, defined by  Ozsv\'{a}th and Szab\'{o} \cite{OS:Knot2003}  and independently Rasmussen \cite{Rasmussen2003}, that comes in several different flavors.  The hat version $\widehat{\HFL}$ is the simplest and $\HFL^-$  contains more geometric information.
$\HFL$ detects the genus \cite{OSGenusBounds} and fiberedness \cite{Ni2007, Ni2009, Ghiggini2008} of knots  and  categorifies the Alexander polynomial.
That is, the graded Euler characteristic of $\HFL(L)$ is the Alexander polynomial of $L$.  
A combinatorial method of calculating link Floer homology was presented in \cite{MOS2006, MOT2009}.  
In \cite{MOST2006}, a completely combinatorial construction of $\HFL$ was given that was independent of holomorphic techniques and used a grid diagram representation of a link as input.

The first evidence that markings were needed came when Juh\'{a}sz \cite{Juhasz2009}, using sutured Floer homology (another variation of $\HFL$), showed that $\widehat{\HFL}$  is functorial with respect to smooth decorated cobordisms.  
One might have  hoped that these markings would only be required when working with sutured Floer homology and/or when dealing with the hat version  $\widehat{\HFL}$.  However, it turned out that markings are important in the combinatorial version as well, which calculates $\HFL^-$.  After defining chain maps on $\HFL^-$, whose underlying grid diagram maps could be viewed as births, deaths and saddles, Sarkar \cite{Sarkar2011} showed that the chain maps induced by his grid diagram maps depended on a particular type of marking.  Indeed, he showed  that rotating this marking once around a knot induces a non-trivial automorphism of the link Floer homology chain complex for most of the 85 prime knots up to nine crossings.  

The strongest evidence for the need of markings comes from the following recent result by Juh\'{a}sz and Thurston \cite{JT2012}.  
Let a \emph{based oriented link} be an oriented  link $L\subset S^3$ along with a set of base points $\mathbf{p}=\set{p_1, \ldots, p_n}$ exactly one on each link component of $L$. Let $\Link_*$ be the category whose objects are based oriented links and whose morphisms are  diffeomorphisms of $S^3$ preserving the based oriented link.
\begin{theorem}
  (Juh\'{a}sz and Thurston 2012) There are functors 
  \begin{align*}
    \widehat{\HFL}&: \Link_*\to \F_2\mbox{-}\Vect,
    \\
    \HFL\phantom{.}\hspace{-.25cm}^-&: \Link_*\to \F_2[U]\mbox{-}\Mod,
  \end{align*}
agreeing up to isomorphism with the link invariants defined by Ozv\'{a}th-Szab\'{o} and Rasmussen.  Isotopic diffeomorphisms induce identical  $\HFL$  maps.
\end{theorem}

For a marked link $L^\bullet=(L,\mathbf{p})$ this means that two different embeddings that are isotopic to one another induce  identical  $\HFL$ maps.  The main reason that link Floer homology has had such great success is the fact that the  isomorphism class of  $\HFL(L)$ is an invariant of the link $L$.  Juh\'{a}sz and Thurston's result strengthened this statement: it shows that $\HFL(L^\bullet)$ is a well-defined group, not just an isomorphism class of groups.

On a separate but related note, $\HFL$ is related to Khovanov homology of a link $L$ via a spectral sequence from the Khovanov homology of $L$, $\Kh(L)$, to the link Floer homology of the branched double cover of $L$ \cite{OSSpectral}.
Khovanov homology is also a powerful link invariant, it: categorifies the Jones polynomial \cite{KhCategorification}; puts a lower bound on the slice genus \cite{Rasmussen2010};  gives an upper bound on the Thurston-Benneqin number of a Legendrian link \cite{Ng2005}; and (a reduced version) detects the unknot \cite{KM2011}.  However, for us, the most interesting property is the fact that Khovanov homology is functorial over smooth isotopy classes of (orientable) surfaces in four space \cite{KhFunctorial}.  Specifically this means: (1) a smoothly embedded (orientable) surface in four space with boundary $L_1\cup L_2$ induces a map $\Kh(L_1)\to \Kh(L_2)$; and (2) two smoothly embedded (orientable)   surfaces in four space that are isotopic induce identical maps $\Kh(L_1)\to \Kh(L_2)$.

Combining the apparent need of markings in link Floer homology with the well defined correspondence between Khovanov homology and link Floer homology along with the fact that Khovanov homology is functorial over smooth isotopy classes of surface in four space, we pose the following conjecture.

\begin{conjecture}
  Link Floer homology is functorial with respect to marked smooth isotopy classes of surfaces embedded in four space.  Specifically this means: (1) a marked surface smoothly embedded in four space with boundary $L_1^\bullet\cup L_2^\bullet$ induces a map $\HFL(L^\bullet_1)\to \HFL(L^\bullet_2)$; and (2) two different marked surfaces embedded in four space that are isotopic induce identical maps $\HFL(L^\bullet_1)\to \HFL(L^\bullet_2)$.
\end{conjecture}

To get the analogous result in Khovanov homology, Khovanov \cite{KhFunctorial} appealed to a movie move theorem, which gave him a list of moves that he needed to verify held certain properties.  Specifically, each movie move induces two different chain maps, one for each side of the movie move.  It needs to be verified, for each movie move, that the two induced chain maps are in fact chain homotopic in order to claim that isotopic embeddings induce the same map on homology.

One would like to use a similar approach for link Floer homology.  
However, to date, there has been little work done on marked isotopy classes of surfaces in 4-space.

The objective of this paper is to find a reasonable definition of marking of surfaces, and isotopy of marked surfaces, that is compatible with the markings needed for link Floer homology.  Specifically, we present a marked analog of the following theorem of Carter and Saito \cite{CS1993}.

\begin{figure}[h]
  \centering
  \includegraphics[width=\textwidth]{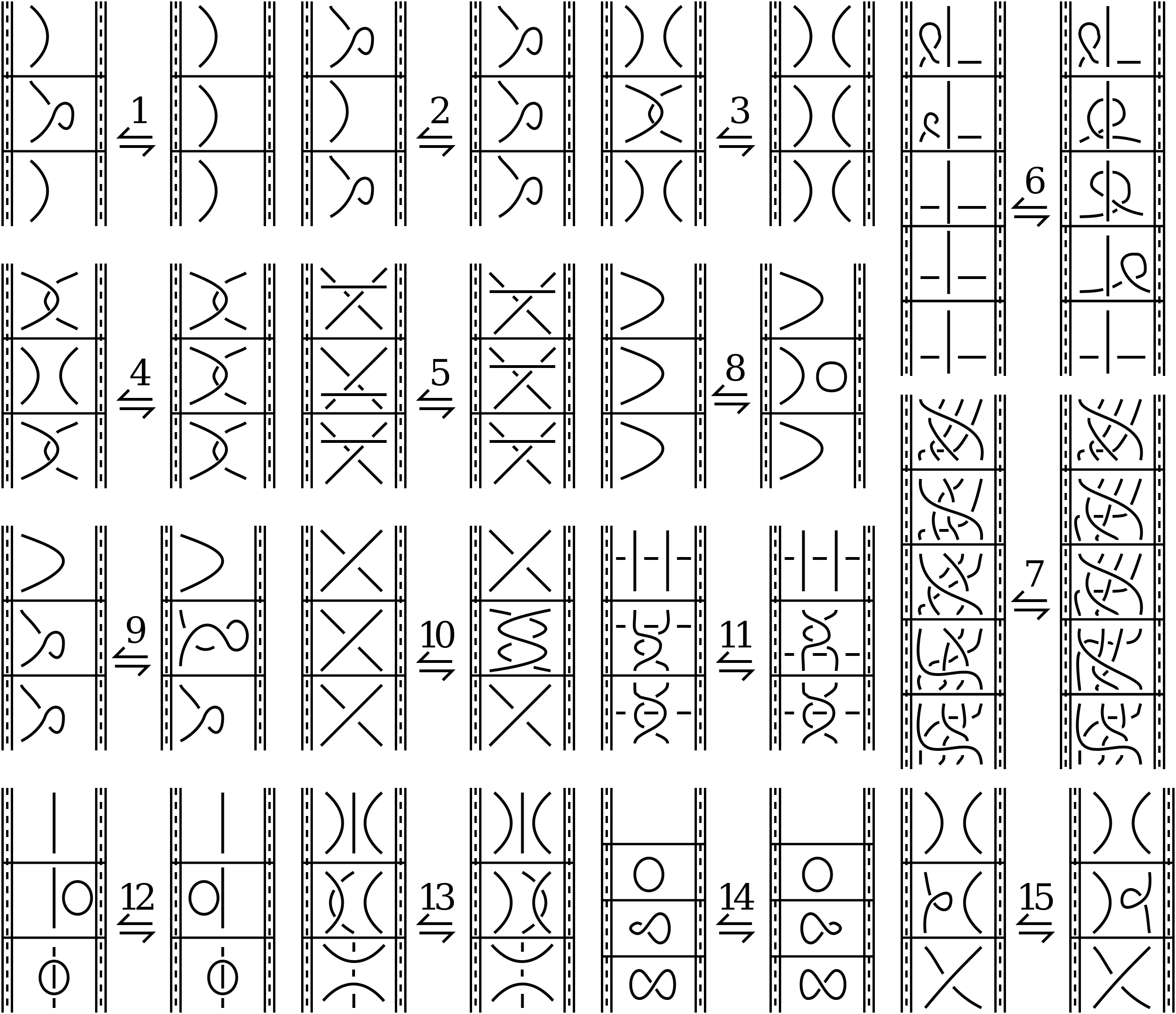}
  \caption{The fifteen CS-movie moves.  }
  \label{fig:SmoothMovieMoves}

\end{figure}

\begin{theorem}
 (Carter and Saito 1993)
  Two knotted surface movies represent isotopic knottings if and only if they are related by a finite sequence of the 15 movie moves (depicted in figure \ref{fig:SmoothMovieMoves})
or interchanging the levels of distant critical points.
   \label{theorem:Smootheoremovie}
\end{theorem}

It should be noted that there are two different movie  theorems.  The CS-movie theorem just stated uses a smooth collection of stills to represent a particular embedding of a surface and uses the 15 movie moves in figure \ref{fig:SmoothMovieMoves}.  The CRS-movie theorem \cite{CRS1997},  requires a second height function in each still.  It uses a finite collection of stills to represent an embedded surface and has 31 different movie moves.

We do not discuss a marked version of the CRS-movie move theorem in this paper although we intend on supplying such a theorem in the future.

The rest of this paper is organized as follows.
In the next section, we give a brief review of movies with boundary.  After this we discuss markings and  define marked surfaces and three marked equivalences that are necessary in the   proof of the  marked analogue of the marked CS-movie theorem that appears in the final section.

\section{Movies with Boundary}
In this section we give a brief review of movies with  boundary   in order to fix notation used in the rest of the paper.  Further details can be found in \cite{CRS1997, Graham2013}.

Let $(F,\partial F)$ be a surface with boundary.  Choose a vector $v_1\in \R^4$ and for each $s\in \R$ let $\R^3_s$ denote the hyperplane orthogonal to $v_1$ at the point $sv_1$.
\begin{definition}
  An embedding of a knotted surface with boundary $(K,\partial K, v_1, a,b)$ is an embedding $K\colon F \hookrightarrow\R^4$  along with a choice of  unit vector $v_1\in \R^4$  and $a<b\in \R$ such that
  \begin{enumerate}
  \item  $\partial K\subset \R^3_a\cup \R^3_b$;
  \item $K\subset \R^3\times[a,b]$;
  \item $\left (K-\partial K\right) \subset \R^3\times (a,b)$.
  \end{enumerate}
\end{definition}
As in knot theory, we will let $K$ refer either to the embedding $K\colon F\hookrightarrow \R^4$ or the image $K(F)$ and $\partial K :=K(\partial F)$.

 \begin{definition}
     Two knotted surfaces with boundary $(K_i, \partial K_i, v_1, a, b), \ i=1,2$  are \emph{boundary ambiently isotopic} if there is an isotopy $H\colon \R^4 \times [0,1] \to \R^4$ such that: 
     \begin{enumerate}
       \item $H(x,0)=x \ \ \forall x\in \R^4$;
       \item $H(K_1(c), 1)=K_2(c) \ \ \forall c\in F$;
       \item $H(K,s)\subset \R^3\times[a,b]  \ \forall s\in [0,1]$;
       \item $\left.H_s\right|_{\R^3_a\cup \R^3_b}= \R^3_a\cup \R^3_b \ \ \forall s\in [0,1]$.
     \end{enumerate}
   \end{definition}

A \emph{knotted surface diagram}, $\K:=p_v(K)$, is the image of $K$, under a projection orthogonal to $v$, to an $\R^3$ hyperplane of $\R^4$ along with  a depiction of crossing information.  The orthogonal projection $p_v\colon \R^4\to \R^3$ can be used to keep track of relative height information of double points, triple points, etc.
The composition $p_v\comp K(F)$ is \emph{generic} if it is an immersion except for some isolated points, which are cone points of figure 8's.
The \emph{$j$-tuple set} is $S_j=\set{y\in R: \# (p_v\comp K)\inv(y)=j}$.  For $j=2,3$ this is the double point and triple point set respectively.  
A \emph{branch point} is a point $y\in \R^3$ such that the intersection of  any neighborhood $N(y)$ with $p_v\comp K(F)$ contains a cone on a figure 8.  
The closure of the double point set $\overline{S_2}$ contains the branch points and the triple points and can be defined as the image of a compact 1-dimensional manifold  (non-generically immersed) in $\R^3$ (for details see \cite{CRS1997, Graham2013}).
From now on we include the branch points in the double point set even though the  preimage of a branch point is a single point.

Projection onto the vector $v_1$, $p_1\colon \R^3\to \R$ is a \emph{generic Morse function for the knotting} (comprised of fixed embedding $K$ and generic map $p\comp K$) if: (1)
the critical points, with respect to $p_1$, of $S_1, S_2, S_3$,   are non-degenerate; and (2) each critical point is at a distinct critical level of $p_1$. 
We define critical points to include triple points and branch points.
Furthermore, critical points of type I, II and III refer to the critical points of $S_1, S_2$ and $S_3$ respectively.  Type I critical points are the standard index 0, 1 and 2 critical points corresponding to births, saddles and deaths respectively.

A  \emph{CS-movie} of the knotted surface $K$ is a knotted surface diagram $\K$ with a fixed generic Morse function.
However, in the case of knotted surfaces with boundary we chose the vector $v_1$ before the vector $v$ and therefore may need to perturb $v_1$  to get the required Morse function $p_1$.  
If a perturbation is needed we  need to perturb the boundary components of $K$ to lie in $\R^3_a$ with respect to $v_1'$ instead $v_1$ to ensure that we are still working with an embedded surface with boundary 

A \emph{complementary coordinate system} for an embedding $K$ is an orthonormal coordinate system of $\R^4$,  $(v,v_1,v_2,v_3)$ that satisfies: (1) the projection orthogonal to  $v$,  $p_v\colon\R^4\to \R^3$ is generic with respect to $K$;  (2) the  projection  of $p_1\colon\R^4\to \R$ onto the vector $v_1$ is a generic Morse function for $K$; and  (3) $p_1$ is a generic Morse function for $\K$.
A complimentary coordinate system using the vector $v_1$ and $a>b\in \R$ can always be chosen for any knotted surface                  (see \cite{Graham2013}).
Furthermore, in \cite{Graham2013}, we defined isotopies of knotted surfaces with boundary that take a knotted surface with boundary with choices $v_1$ and $a>b\in \R$ to some other embedding with $v_1'$ and $a'>b'\in \R$ .
Therefore, for simplicity, we will always assume when considering  multiple marked smooth movies that a single complimentary coordinate system is chosen with some fixed $v_1$ and $a>b\in \R$.

\section{Equivalence of Markings}

In this section we give a definition of a marking of an embedded surface in $\R^4$ and then discuss the complications that arise when isotopies of the surface are considered.  We then circumvent these complications by defining three equivalence relations for markings and work over an equivalence class of marked surfaces that will be used in the following section to prove a marked movie move theorem. 

\begin{definition}
  A \emph{marking} of an embedded surface $K\subset \R^3\times[a,b]$ is an embedding of a connected graph $T$ with only univalent and trivalent vertices $\xymatrix{T \incl[r]^M & K \subset \R^4}$ such that:
  \begin{enumerate}
    \item There exists a single marked point in each link component of $K\cap \R^3_s$ for each $s\in [a,b]-\set{\mbox{ type I critical values }}$.
    \item 
      \begin{enumerate}
      \item If $x$ is a single valence vertex then  $p\left( M(x)\right)$ is either an index 0 or 2 type I critical point or a marking of a link component in $p(\partial K)$;
      \item If $x$ is a trivalent vertex then  $p\left( M(x)\right)$   is an index 1 type I critical point;
      \end{enumerate}
  \end{enumerate}
\end{definition}
A \emph{marked surface} is an embedded surface with a marking.  Again, we will let $M$ refer either to the embedding or the image of the embedding $M(T)$.
\begin{definition}
  A \emph{flow orientation} of a marking is the orientation that is endowed by flowing in the positive time direction along the marking.
\end{definition}
Now that we have marked surfaces we need to define what a marked isotopy of these surfaces is.  
Suppose that we let a marked isotopy be a level preserving isotopy $g$ of $K$.  
Then $g$ automatically is an isotopy through markings since each link component in each level set has a single marking and the level sets are preserved.  
However, level preserving isotopies are quite a restrictive class, especially in light of the fact that the  non-marked marked CS-movie theorem is a theorem about arbitrary isotopies.

\begin{figure}[h]
  \centering
  \includegraphics[width=\textwidth]{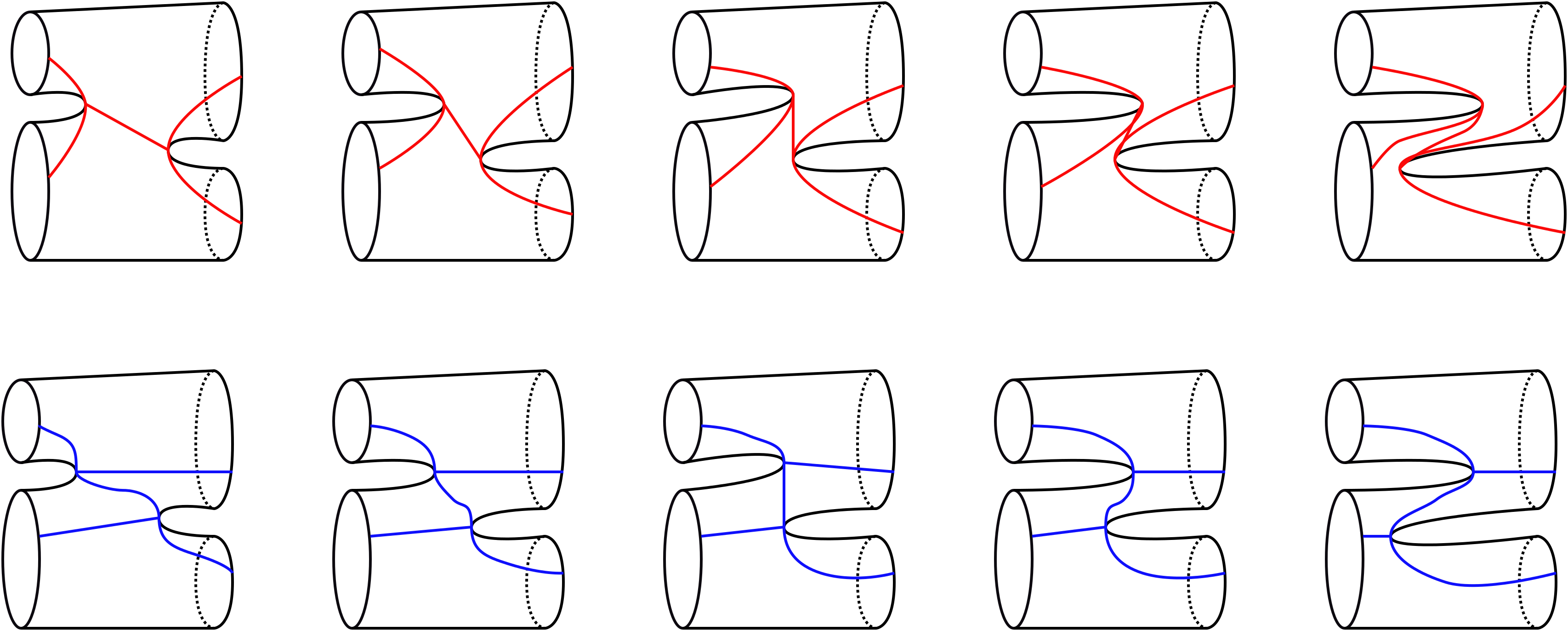}
  \caption[ \ \ Isotopy exchanging a merge and split saddle critical point]{An isotopy exchanging a merge and split saddle point.  The top left two diagrams (and the bottom right two diagrams) depict an isotopy acting on a marked surface.  The middle diagrams depict when the two critical points lie in the same level set.  The sets of marked points in the top right three (and the bottom left three) diagrams are not markings since there are more than one marking per link component in some level sets.}
  \label{fig:ProblemIsotopy}
\end{figure}

Alternatively, we could require a marked isotopy to be an isotopy $g$ of $(K,M)$ such that $(g_s(K), g_s(M))$ is a marked surface for all $s$ (\ie it is an isotopy through markings).
At first glance, this seems like a good definition since it allows the critical values to change.
However, there are some type I critical values (the Morse critical points) that must maintain their ordering.
For example, in many cases exchanging the ordering of two type I index 1 critical values cannot be done through markings.
A specific example of this is in figure \ref{fig:ProblemIsotopy}.  

Consider the top sequence of diagrams in figure \ref{fig:ProblemIsotopy} that represent an isotopy $g$ acting on a marked surface $(K,M)$ and assume that the middle  diagram occurs at $g_\h $ and  $g_0, g_1$ are the  left most and right most  diagrams respectively.  
Further assume that $x_i$ and $x_j$ are the two type I index 1 critical points depicted, with critical values $c_i<c_j$ in $g_0$ (according to our  choice of vector $v_1$). 
Then $g_s(K,M)$ is a marked surface for $s\in [0,\h)$.  
However, $g_{\h}(K,M)$ is not a marked surface since $g_{\h}(M)$ has an arc of the marking $\varphi$  residing in some  link component in the level set $c_i=c_j$.  
The image of the marking in $g_s(M)$ for $ s\in (\h,1]$ is not a marking since in each still between the two critical values $c_j<c_i$ there is a link component with exactly three markings.

Define $\og_s=g_{1-s}$ (\ie $\og$ is the isotopy that runs $g$  backwards), which means $\og_0$ is the right most diagram of the second sequence.  Interestingly, there is another marking $M'$ of $g_1(K)$ pictured in the second sequence that has the exact same properties as $M$.  That is to say  $\overline{g}(g_1(K),M')$ has the exact same properties as $g(K,M)$.
Namely,
\begin{enumerate}
  \item  $\og_s(g_1(K),M')$ is a marking for $s\in [0,\h)$; 
  \item there is an arc $\varphi'$ in $\og_{\h}(M')$ residing in  some link component in the level set $c_i=c_j$ (as oriented arcs $\og_{\h}(\varphi')=-g_{\h}(\varphi)$, that is they are the same arc with opposite orientation (assuming that $g_0(M)$ and $\overline{g}(M')$ are given flow orientations--left diagram of figure \ref{fig:MarkingExchange});
  \item the image of the marking $\og_s(M')$  is not a marking since in each still between the two critical values $c_j<c_i$ there is a link component with exactly three markings.
\end{enumerate}

It is clear that this definition of marking does not always allow for isotopy through markings when exchanging critical point levels.
Exchanging critical point levels is essential to the smooth unmarked movie move theorem.
This means that we need to find a notion of marking and marked isotopy that can always exchange critical points in a well-defined manner.
To do this we will modify the class of markings by defining equivalence relations.
The first equivalence relation defined in the next subsection equates $g_0(K,M)$ with $\og_1(g_1(K),M')$.

\subsection{Marking Equivalences}
The point of this subsection is to modify the notion of  equivalence for  markings so that there is a sensible way to exchange the levels of critical points  while staying within the same equivalence class of marking.
We will define three different equivalence relations: a (balancing) twist equivalence relation; a canceling equivalence; and a zero arc equivalence.
The arcs $\varphi, \varphi'$ encountered in the example depicted in figure \ref{fig:ProblemIsotopy} that connect the two type I index one critical points are special arcs that we will name zero arcs.

\begin{definition}
  Let $x_i$ and $x_{j}$ be two successive type I index 1 critical points,  with critical values $c_i<c_j$, and let $g\colon \R^4\times [0,1]\to \R^4$ be an isotopy  of a marked knotted surface $(K,M)$ such that $g_1(x_i)$ and $g_1(x_{j})$ are contained in the same level set and $g_s(M)$ is a marking for $s\in[0,1)$.  
A \emph{zero arc} or \emph{zero curve}, of the marking $M$, $z_{ij}\subset g_0(M)$ connects the critical point $x_i$ to  $x_{j}$ and is an arc such that $g_1(z_{ij})$ is an embedded arc, of the link with two singular crossings, that intersects  $g_1(x_i)$ and $g_1(x_{j})$ at its endpoints and nowhere else.
\end{definition}

\begin{figure}[h]
  \centering
  \includegraphics[width=.45\textwidth]{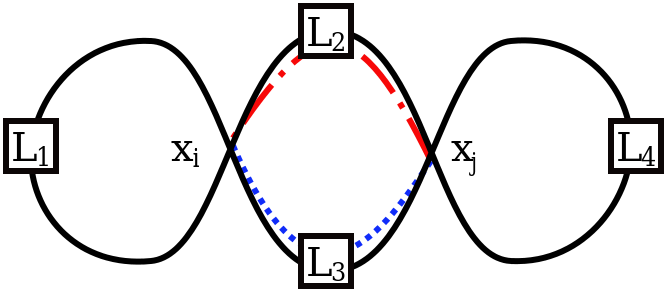}
  \caption[ \ \ Zero arcs]{The image of two zero arcs, represented by dotted lines, connect two singular points $c_i$ and $c_{j}$.  $L_i$ are  different  links.}
  \label{fig:ZeroArcs}
\end{figure}

There  are only  two possible images of  zero arcs and they  are depicted in figure \ref{fig:ZeroArcs}.
Suppose that initially  $c_i<c_{j}$,  then the zero curve $z_{ij}$ would be oriented with initial point $x_i$ and final point $x_{j}$ and likewise $g_1(z_{ij})$ has initial point $g_1(x_i)$ and endpoint $g_1(x_{j})$.
If $ c_{j}<c_i$ initially the orientations are  reversed and we would denote this by the zero arc $z_{ji}$.
There are four different cases for passing type I index 1 critical points past each other, depending on whether $x_i, x_{j}$ are merge or  split saddle points.  
We have discussed an example involving a merge and split type I index 1 critical points depicted in figure \ref{fig:ProblemIsotopy}.

\subsubsection{Zero Arc Equivalence}
\begin{figure}[p]
  \centering
  \includegraphics[width=\textwidth]{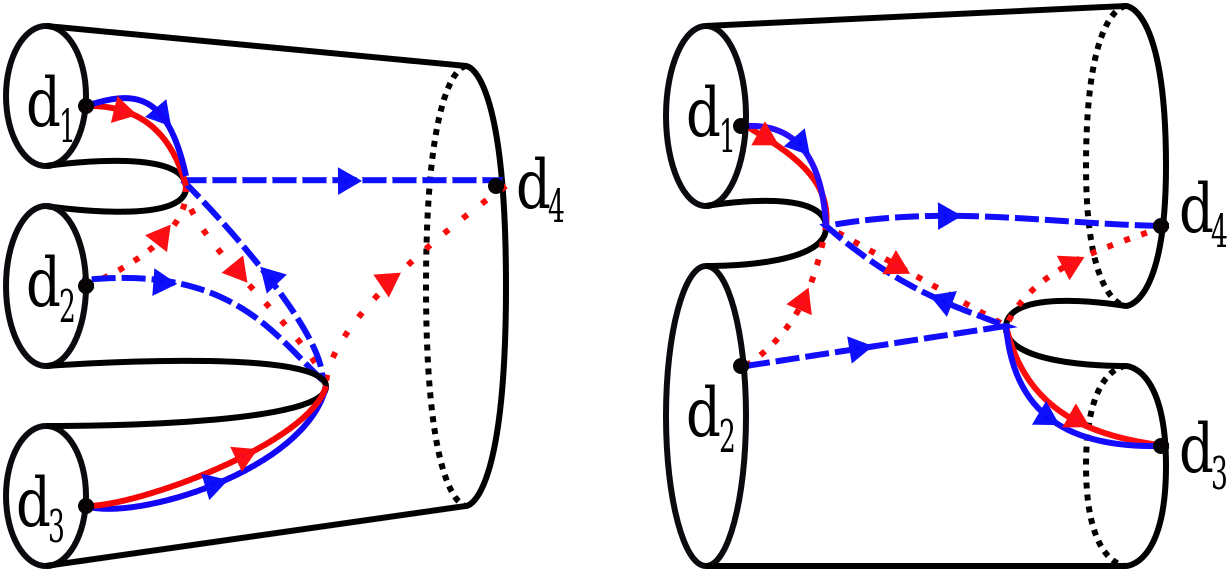}
  \caption[ \ \ Zero arc equivalence 1]{Diagram of curves used in the definition of zero arc equivalence.  Dotted arcs do not represent arcs going around the back (even though in general they could do so). They are curves on the front of the diagram.}
  \label{fig:MarkingInitial}
\vspace{.5in}
  \centering
  \includegraphics[width= \textwidth]{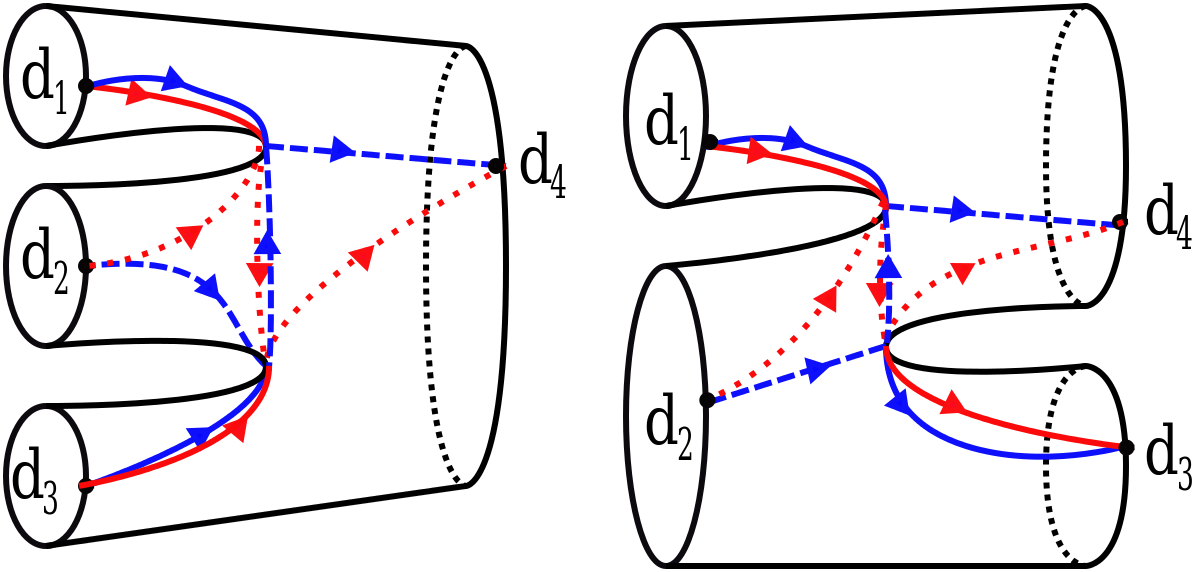}
  \caption[ \ \ Zero arc equivalence 2]{Two different examples of the moment in the isotopy $h$ that has both critical values in the same level set.  The light (red) marking defined in $\lim_{s\to \h}h|_{[0,s)}$  is equivalent to the dark (blue) marking defined in $\lim_{s\to \h}h|_{(s,1]}$ through the isotopy $h$ in each case.  Again the dotted lines do not represent arcs going around the back.}
  \label{fig:MarkingExchange}
\end{figure}

We would like to say that the markings in the top two diagrams of figure \ref{fig:ProblemIsotopy} are equivalent to the markings in the bottom right two diagrams.  
Then if we are dealing with an isotopy that exchanges index 1 critical points we can stay within the same equivalence class of marking except at the point in the isotopy in which both critical points lie in the same level set.  
The following constructions make this equivalence precise.
Follow figures \ref{fig:MarkingInitial} and \ref{fig:MarkingExchange} when making sense of the following definition.

\begin{definition}
  Let $g\colon \R^4\times[0,1]\to\R^4$ be an isotopy that exchanges two consecutive type I index 1 critical points $x_i,x_j$, with critical values $c_i, c_j$, of a boundary knotted surface  $K: =K_{[a,b]}$ and let $\og_s: =g_{1-s}$ be the isotopy that runs $g$ backwards.  Let $M$ be a marking of $g_0(K)$ and  $M'$ be a marking of $g_1(K)$.  Then the markings $M$ and $M'$ are \emph{zero equivalent through $g$} if they satisfy the following (suppose $x_i,x_j$ only share a level set at $g_\h$ and $\epsilon>0$ is small):

  \begin{enumerate}
    \item $g|_{s}(M)$ and $\og|_s(M')$ are markings for $s\in[0,\h)$;
    \item $M$ has a zero curve $z: =z_{ij}\subset g_0(M)$ connecting  $x_i$ to $x_j$;
    \item $M'$ has a zero curve $z': =z_{ji}\subset \og_0(M')=g_1(M')$ connecting  $x_j$ to $x_i$;
    \item  $g_\h(z)=-\og_\h(z')$ (\ie they are the same arc with opposite flow orientations);
    \item there are four points in the  intersection $M\cap \og_1(M')$ residing in $\R^3_{c_i-\epsilon}\cup\R^3_{c_j+\epsilon}$
labeled by $d_k$ for $k=1,\ldots, 4$ with at least one intersection point in each level set $c_i-\epsilon$ and $c_j+\epsilon$;
    \item $g_s\left(M_{[a,c_i-\epsilon]}\right)$ is a marking for $g_s(K_{[a,c_i-\epsilon]})$ and  $g_s\left(M_{[c_j+\epsilon,b]}\right)$ is a marking for $g_s\left(K_{[c_j+\epsilon, b]}\right)$ for each $s \in [0,1]$ and $g_1\left(M_{[a,c_i-\epsilon]}\right)=M'_{[a,c_i-\epsilon]}$  and $g_1\left(M_{[c_j+\epsilon,b]}\right)=M'_{[c_j+\epsilon,b]}$.
    \item each $d_k$ is an endpoint of two arcs $e_k\subset M_{[c_i-\epsilon,c_j+\epsilon]}$ and $e_k'\subset \overline{h}_1\left(M'_{[c_i-\epsilon,c_j+\epsilon]}\right)$ each containing only one of $x_i,x_j$, which is the other endpoint;
      for exactly two values of $k=1,\ldots, 4, \ e_k$ is isotopic to $e_k'$ relative to the endpoints (in figure \ref{fig:MarkingInitial} these arcs are solid and have endpoints in $d_1$ and $d_3$);
    \item for the other two values, say $l,m$, of $k$ not used in the last condition there are two  arcs with endpoints $d_l$ and $d_m$ which have both $x_i$ and $x_j$ in the interior;  
      these two arcs are isotopic (these are the dotted arcs with endpoints in $d_2$ and $d_4$ in figure \ref{fig:MarkingInitial}).
  \end{enumerate}
\end{definition}

Condition 5 allows us to separate the marked surface into different regions: an interesting middle region that is important for the equivalence relation; and  exterior regions.
Condition 6 states that the markings outside the interesting region are isotopic to one another.
Condition 7 and 8 guarantee that the curves we would like to make equivalent  inside the level sets $[c_i-\epsilon, c_j+\epsilon]$ are isotopic.

\begin{figure}[h]
  \centering
  \includegraphics[width=\textwidth]{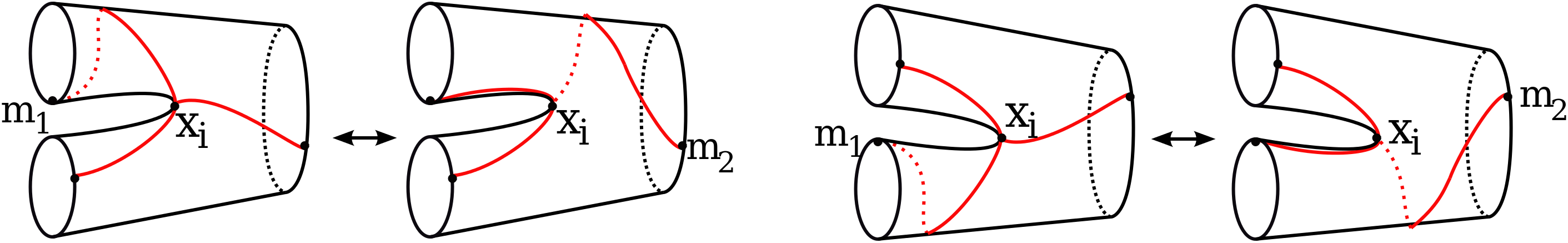}
  \caption[ \ \ Twist equivalent markings]{Examples of twist equivalent markings the marked points $m_1,m_2$ and $x_i$ are defined in  \ref{definition:TwistEquivalence}.  The isotopy is given by flow in the positive time direction.}
  \label{fig:TwistEquivalentMarkings}
\end{figure}

\FloatBarrier
The needed zero arcs for an isotopy $g$ of a marked surface $(K,M)$ that exchanges index 1 points may not exist.  
This means there is not a zero equivalent marking to $M$ through $g$.
This simply means that the arc connecting the two critical points does not have the appropriate number of twists. The next equivalence relation will  allow us to move twists past critical points in order to create the needed zero arc.  See figure \ref{fig:TwistEquivalentMarkings} for the motivating picture of the twist equivalence.

\subsubsection{Twist Equivalence}
Let $(K_{[a,b]},M)$ be a marked surface with indexed critical levels $c_0<c_1<\cdots<c_n$.  Suppose that $x_i$ is a type I index 1 critical point with critical value $c_i$.  Let $m$ be a point of the marking in the $a_i$ level where $c_{i-1}<a_i<c_i$.  Let $w_i\subset M\cap K_{[a_i,c_i]}$ be  the marked  subarc of $M$, with endpoints $m$ and $x_i$, 
\begin{definition}
  Let $g\colon \R^4\times[0,1)\to \R^4$ be an isotopy of the marked surface $(K_{[a,b]},M)$, defined above, such that: $g_s(K)=K$ for each $s\in[0,1)$; $g_0(m)=m\neq x_i, \ \lim_{s\to 1}g_s(m)=x_i$ and $g$ is the identity when restricted to the level sets ${[a,b]-[a_i,c_i)}$.  Then the arc $w_i: =M\cap K_{[a_i,c_i]}$ is a (single) \emph{twist arc of $x_i$} with respect to $g$ if  $g_s(M)$ is a marking for each $s\in[0,1)$ and  $\lim_{s\to 1}g_s(w_i)$ is an immersed oriented arc, in the critical level $c_i$, with exactly one double point $x_i$.
\end{definition}
\begin{definition}
\label{definition:TwistEquivalence}
  Let $K_{[a,b]}$ be a knotted surface with a type I index 1 critical point $x_i$, with critical value $c_i$.  Let $M^1$ and $M^2$ be two markings each containing a single twist arc $w_1 \subset M^1\cap K_{[a_1,c_i]}$  with respect to $g\colon \R^4\times [0,1)\to \R^4$ and $w_2\subset M^2\cap K_{[c_i,a_2]}$ with respect to $h\colon\R^4\times [0,1)\to \R^4$.
$M^1$ and $M^2$ are \emph{twist equivalent through $g$ and $h$} if there exists $a_1,a_2 \in \R$  with  $a<a_1<c_i<a_2<b$ 
such that:

\begin{enumerate}
  \item $M^1_{[a,a_1]}$ is isotopic to   $M^2_{[a,a_1]}$ and $M^1_{[a_2,b]}$ is isotopic to  $M^2_{[a_2,b]}$;
  \item $\lim_{s\to 1}g_s(m_1)\cap \R^2_{c_i}=\lim_{t\to 1}h_t(m_2)\cap \R^2_{c_i}=x_i$ 

(where  $m_1=w_1\cap \R^2_{a_1}$ and $m_2=w_2\cap \R^2_{a_2}$) 
  \item  $\lim_{s\to 1}g_s(w_1)=\lim_{t\to 1}h_t(w_2)\subset \R^2_{c_i}$ as oriented arcs;

\end{enumerate}
\end{definition}

\begin{figure}[h]
  \centering
  \includegraphics[width=.8\textwidth]{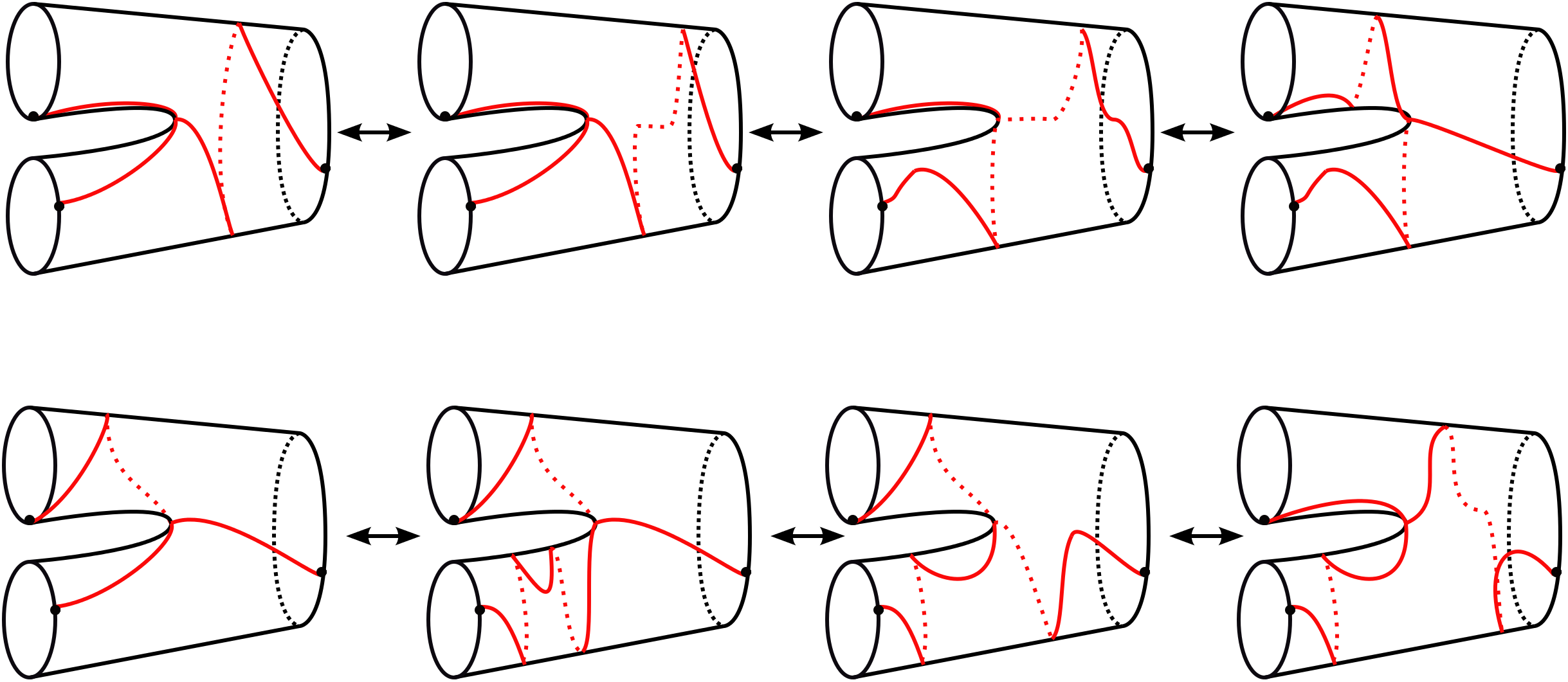}
  \caption[ \ \ Moving twists past saddle critical points]{Two sequences of twist equivalences moving twists past a saddle point.  The top sequence moves the twist right to left assuming an isotopy that flows  in the negative time direction.  The bottom sequence moves the twist left to right assuming an isotopy that flows linearly in the positive time direction (note that this cannot be done in one step since the marking enters and exits the critical point from different directions, which makes the third condition of twist equivalence fail).}
  \label{fig:TwistEquivalenceSequenceExamples}
\end{figure}
\subsubsection{Canceling Equivalence}
An isotopy that cancels either a zero or two index critical point with an index 1 critical point also creates problems for markings.
Consider an isotopy $g\colon \R^4\times[0,1]\to \R^4$ of a marked surface $(K,M)$ that cancels an index zero and one critical point and a marking $M'$ for $g_1(K)$.
There is no hope that $M$ and $M'$ will be isotopic since $M$ has one more univalent and trivalent vertex than $M'$.
The following equivalence relation deals with these cases.

\begin{definition}
  Two marked surfaces $(K_i,M_i), i=1,2$  are \emph{canceling equivalent} through an isotopy  $g\colon \R^4\times[0,1]\to \R^4$ if the following are satisfied:

  \begin{enumerate}
    \item at $s=1$, $g$ cancels an index zero (or two) critical point $x_i$ with an index one critical point $x_j$ with initial critical levels $c_i<c_j$ (correspondingly $c_i>c_j$);
    \item $g_0(K_1)=K_1$ and $g_1(K_1)=K_2$;
    \item $g_s(M_1)$ is a marking for $g_s(K_1)$ for $s\in[0,1)$;
    \item there exists an arc $b_{ij}$ in $g_1(M_1)$ lying in the critical level $g_1(c_i)=g_1(c_j)$ with endpoints $g_1(x_i)$ and $g_1(x_j)$;
    \item the closure of $(g_1(M_1)-b_{ij})$ is isotopic to $M_2$.
  \end{enumerate}
\end{definition}

\begin{figure}[h]
  \centering
  \label{fig:CancelingEquivalence}
  \includegraphics[width=\textwidth]{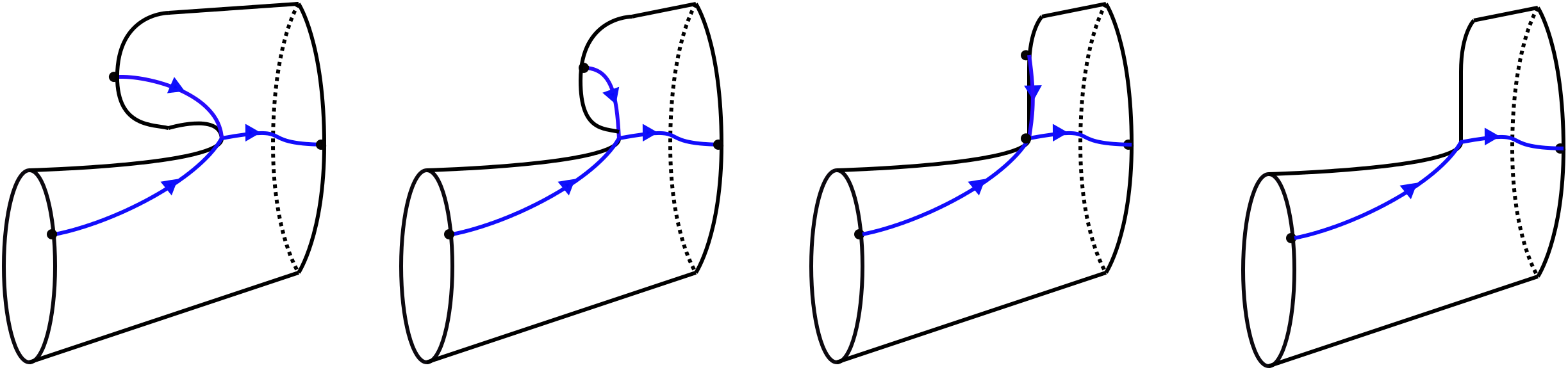}
  \caption[ \ \ Canceling equivalence example]{An example of a canceling equivalence.  The initial marked surface $(K_1,M_1)$ is canceling equivalent to the final marked surface $(K_2,M_2)$ through an isotopy that cancels the index 0 and 1 type I critical points.}
\end{figure}

\begin{definition}
  Let $(K,\overline{M})$ be the equivalence class of markings for $K$  where  $\overline{M}$ is the equivalence class of markings generated by isotopy through markings  and twist, canceling and zero arc equivalences.
\end{definition}

\begin{definition}
  Two marked knotted surfaces $(K_i,M_i)  
,\ \ i=1,2$ are  \emph{marked isotopic} if there exists a smooth isotopy  $g\colon \R^4\times [0,1] \to \R^4$ such that the following three conditions are satisfied.
  \begin{enumerate}
    \item  $g_0(K_1)=K_1$ and $g_1(K_1)=K_2$    
    \item  $(g_1(K_1), g_1(M_1))$ is in the same equivalence class as $(K_2,M_2)$
    \item  if $\partial K\neq \nullset$ and $\partial K\subset \R^3_a\cup \R^3_b$ then $g(\partial K,s)=\partial K$  for each $s\in [0,1]$.
  \end{enumerate}
\end{definition}

\section{Marked Smooth Movie  Theorem}
Now that we have defined the equivalence classes of markings for knotted surfaces, we are in a position to generalize Carter and Saito's movie move theorem to the category of  marked smooth surfaces $K$ with markings in $\overline{M}$.  The projection of a marked knotted surface  $(K,M)$ by  generic function $p\colon \R^4 \to \R^3$ onto the hyperplane $\R^3$, along with crossing information, is a  \emph{marked knotted surface diagram}, denoted  by $(\K,\M): =(p(K),p(\overline{M}))$. 
A \emph{marked knotted surface movie} or simply \emph{marked CS-movie} is a marked knotted surface diagram with a choice of Morse function $p_1\colon \R^3\to \R$ that separates the type I, II, III  critical points except the endpoints of the double point set in $\partial K$.

We denote the unmarked CS-movie moves that operate on  a movie ${\mathcal K}$ by $\MM_i$ for $i=1, \ldots, 15$.  The following lemmas are used to  define marked  movie moves $\MM^\bullet_i$.

\begin{lemma}
\label{lemma:LevelPreservingFlow}
  Suppose that $\K_{[a_1,b_1]} \subset (\K,\M)$ contains no type I singularities and suppose that  $g_s$ is an isotopy through markings satisfying $g_t(\K)=\K$ for each $t\in [0,1]$ defined by flowing the marking uniformly in the time direction.  Then there exists a level preserving isotopy $g$ such that $g_t(\K)=\K$ for all $t\in [0,1]$ and $g_t(M_{[a_1,b_1]})=g_t(M_{[a_1,b_1]})$ for each $t\in [0,1]$.
\end{lemma}
\proof
Let $\K_{[a_1,b_1]}$ and $g_t$ have the properties stated  above. 
Define an isotopy, $P^c_t:\K_c\to \K_c$, to be the process of flowing all points in  the level set $\K_c$  in a direction  orthogonal to $v$, by $P^c_t(\M_c):=g_t(\M_{[a_1,b_1]})\cap \K_c$.
The fact that $\M$ is a marking, there are no type I critical points in $\K_{[a_1,b_1]}$ and that $g_t$ is an isotopy through markings guarantee that $g_t(\M_{[a_1,b_1]})\cap \K_c$ always consists of a single point, which defines the orthogonal flow for the rest of the points in $\K_c$.

Next define $g_t:\K_{[a_1,b_1]} \to \K_{[a_1,b_1]}$ by $g_t:=\cup_{c\in[a_1,b_1]} P^c_t(\M_c)$.
By construction $g_t(M_{[a_1,b_1]})=g_t(M_{[a_1,b_1]})$ for each $t\in [0,1]$.
The fact that $g_t$ is an isotopy of $(\K,\M)$ follows from the uniform flow of each level set and the fact that $g_t$ is an isotopy of the marking.

\begin{lemma}
\label{lemma:UnmarkedMovieMoves}
Whenever an unmarked CS-movie move, that does not contain a type I critical point ($\MM=\MM_i$ for $i=1,\ldots,7, 9, 10, 11$), is performed on a marked CS-movie of a marked knotted surface diagram $(\K,\M)$ it can be arranged, using level preserving isotopies and exchanging critical points, so that there are no marked points in the portions of the stills affected by $\MM$.
\end{lemma}
\proof

\begin{figure}[h]
  \centering
  \includegraphics[width=\textwidth]{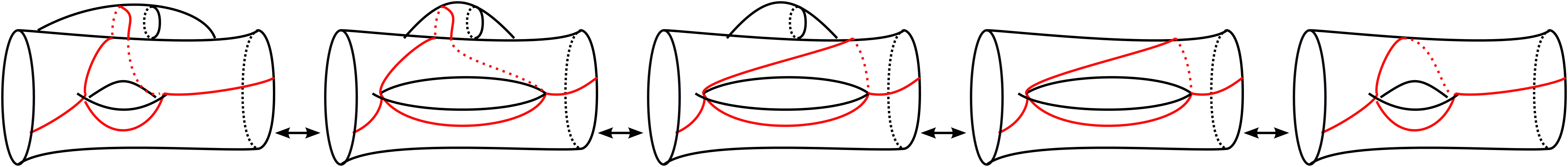}
  \caption[ \ \ Unmarked CS-movie move 1 acting on surface]{Example of unmarked CS-movie move 1 operating on a marked surface.}
  \label{fig:UnmarkedMovieMoveExample}
\end{figure}

Let $(\K,\M)$ and $\MM$ be given as in the lemma with generic projection and Morse function $p_1$.
What one would like to do is to perform a level preserving isotopy that isotopes  the marking away from the unmarked CS-movie move.  
However, this cannot always be done directly.
Sometimes we will need to exchange critical levels of certain critical points.
See the diagram on the left of figure \ref{fig:UnmarkedMovieMoveExample} and note that the type I critical points with values in between the type II branch points obstruct an isotopy through markings.  
The following procedure works in  general.

Let $c_1< \cdots<c_n$ be the $n$ type I,II and III critical points of $p_1\comp p(\K)$. 
Let $\D_{[a_1,b_1]}$ be a collection  of sub-diagrams operated on by $\MM$.
If $\M\cap \D_{[a_1,b_1]}=\nullset$ then there is nothing to prove.
Otherwise, exchange critical points until there are no type I critical points in  $\K_{[a_1,b_1]}$.
This can always be done for this list of movie moves since they do not involve type I critical points.
If there are no type I critical points in $\K_{[a,b]}$ then there exists an isotopy $h\colon \R^4\times[0,1]\to \R^4$ such that $h_s(\overline{M})$ is a marking for each $s\in [0,1]$ and $p(h_1(\overline{M}))\cap \D_{[a_1,b_1]}=\nullset$.

To define $h_s$ let  $\phi\colon \R^4\to \R^4$ be a bump function with support $K_{[a_1-\epsilon,b_1+\epsilon]}$ (for small $\epsilon>0$) and  define $f^1\colon \R^4\times [0,1]\to \R^4$ to be an isotopy that flows uniformly in the time direction satisfying $f_u^1(K)=K$ for  each $u\in [0,1]$ ($f^1_u$ takes diagram 2 to diagram 3 in figure \ref{fig:UnmarkedMovieMoveExample}). 
Since there are no type I critical values in this region  $f^1_u(\M)$ is a marking for $u\in [0,1]$ and by lemma \ref{lemma:LevelPreservingFlow} we know that there is a level preserving isotopy that effects the marking in the same way as $f^1_u(\M)$.
For the rest of the proof assume that we had chosen $f^1_u(\M)$ to be level preserving from the beginning.
Then $p\left(f^1_1\comp \phi( \overline{M})\right)\cap \D_{[a_1,b_1]}$ is either empty or a single connected line segment running in the time direction.  
If it is empty then $h: =f^1\comp \phi$ satisfies the condition.
If it is not empty, then define  $h: =f^2\comp f^1\comp \phi$, where the isotopy $f^2\colon \R^4\times [0,1]\to \R^4$ satisfies $f^2_s(K)=K$ for each $v\in [0,1]$ and flows the points of each level set  $K_t$ in the positive direction (as given by some local orientation with support $K_{[a_1-\epsilon, b_1+\epsilon]}$).
Since $f^2$ is level preserving  the smooth isotopy $f^2_v(f^1_u(\overline{M}))$ is automatically a marking for each $u,v\in [0,1]$.
Then $p\left(f^2_1\comp f^1_1\comp \phi(\overline{M})\right)\cap  B$ is empty. 
Thus in all cases we are able to isotope the markings away from movie moves 1-7, 9, 10, and 11.

In light of this lemma we make the following definition.
\begin{definition}
  The marked CS-movie moves not involving type I critical points are the same as the unmarked CS-movie moves.  That is, $\MM^\bullet_i: =\MM_i$ for $i=1,\ldots,7,9,10,11$.
\end{definition}

This leaves only the movie moves that involve type I critical points.  
However, the definition of marking restricts the possibilities since  each type I critical point is a vertex of the marking.

\begin{figure}[h]
  \centering
  \includegraphics[width=.95\textwidth]{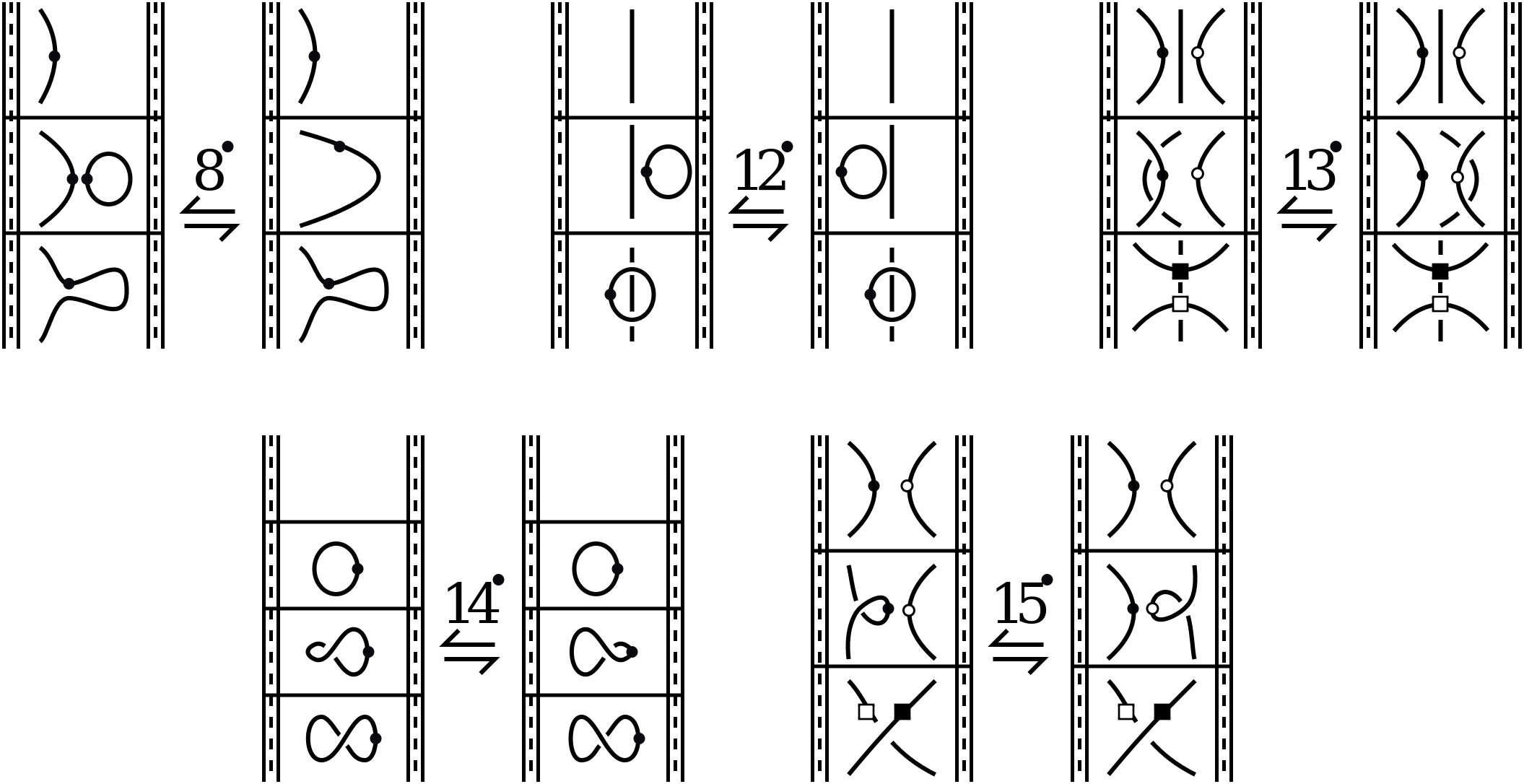}
  \caption[ \ \ marked CS-movie moves]{marked CS-movie moves 8, 12-15. There are various versions of 13 and 15.  Choose either both types of marked circles or both types of marked squares, then choose one of the other type.  This encodes both fission and merge saddle moves.}
  \label{fig:MarkedMovieMoves}
\end{figure}

\begin{definition}
  Marked CS-movie moves 8,12-15 are depicted  in figure \ref{fig:MarkedMovieMoves}.  
\end{definition}

\begin{theorem}
\label{theorem:MarkedMovies}
    Two marked knotted surface movies  $(\K^i, \M^i)$  (for $\K^i: =\K^i_{[a,b]}$ and $i=1,2$) represent   marked isotopic knottings $(K^i,\overline{M}^i)$ if and only if they are related by a finite sequence of marked CS-movie moves $\MM^\bullet_j, j=1,\ldots, 15$, interchanging the levels of distant critical points, twist and zero arc equivalences.
\end{theorem}
\vspace{-.1in}
\begin{proof}
  Suppose that $(\K^i, \M^i)$  represent marked isotopic knottings $(K^i,\overline{M}^i)$ and choose a complementary coordinate system.
As unmarked surfaces,  $\K^i$ represent isotopic knottings $K^i$, which means that $\K^1$ and $\K^2$ are related by a finite sequence of unmarked CS-movie moves $\MM_{l_1}\MM_{l_2}\cdots \MM_{l_n}$ and interchanging levels of distant critical points, where $l_k\in\set{1,\ldots, 15}$ for each $k$.
Perform a level preserving isotopy of the movies so that the markings for $\MM^\bullet_8, \MM^\bullet_{12}, \MM^\bullet_{13},$ $\MM^\bullet_{14}, \MM^\bullet_{15}$ appear as they do in figure \ref{fig:MarkedMovieMoves}.  Also perform any needed level preserving isotopies and exchange of critical points (as in lemma \ref{lemma:UnmarkedMovieMoves}) to make all movie moves $\MM_i^\bullet$ for $i=1, \ldots, 7,9,\ldots,11$ not have any markings on the effected portions of the stills.  Then there exists  the same sequence of marked CS-movie moves  $\MM^\bullet_{l_1}\MM^\bullet_{l_2}\cdots \MM^\bullet_{l_n}$ and the exact same interchanging levels of distant critical points that relate $(\K^1,\M^1)$ to $(\K^2,\M^2)$.
Recall we may need to choose an equivalent representative of the marking when interchanging critical points. 
Additionally, for each canceling equivalence  there is a corresponding marked CS-movie move $\MM^\bullet_8$.

For the converse, suppose that $(\K^1,\M^1)$ and  $(\K^2,\M^2)$ are related by a finite sequence of marked CS-movie moves  $\MM^\bullet_{l_1}\MM^\bullet_{l_2}\cdots \MM^\bullet_{l_n}$, twist and zero arc equivalences and interchanging levels of distant critical points, where $l_k\in\set{1,\ldots, 15}$ for each $k$.  Viewing these as unmarked CS-movie moves  immediately gives that $\K^1$ and $\K^2$ represent (unmarked) isotopic knottings $K^1$ and $K^2$. 

It only remains to show that $\overline{M}^1\equiv \overline{M}^2$.  
However, it is given that $\M^1$ and $\M^2$ are related by twist and zero arc equivalences and  $\MM^\bullet_8$  moves which means that $\M^1\equiv\M^2$ and therefore $\overline{M}^1\equiv \overline{M}^2$.

\end{proof}

\vspace{-.1in}

\subsection*{Acknowledgments}
I am indebted to Daniel Ruberman for the many conversations that eventually lead to the definitions of the equivalences found in this paper.  I also thank Sucharit Sarkar for sharing  the idea of twist equivalence with me.

\bibliographystyle{amsplain}
\bibliography{Research}

\end{document}